\documentclass[12pt]{article}
\usepackage{amssymb}
\usepackage{amsmath}
\usepackage{latexsym}

\begin{document}

\title{Lattice subalgebras of strongly regular vertex operator algebras}
\author{Geoffrey Mason\\
Department of Mathematics \\
UC Santa Cruz}
\date{}
\maketitle

\begin{abstract}
\noindent
 We prove a sharpened version of a
 conjecture of Dong-Mason about lattice subalgebras of a strongly regular
 vertex operator algebra $V$, and give some applications. These  include
  the existence of a  canonical
conformal subVOA $W\otimes G\otimes Z \subseteq V$, and 
a generalization of the theory of minimal models.
\end{abstract}

\bigskip
\noindent
1. Statement of main results\\
2. Background\\
3. Invariant bilinear form\\
4. The Lie algebra on $V_1$\\
5. Automorphisms\\
6. Projective action of Aut$(V)$ on modules\\
7. Complete reducibility of the $V_1$-action\\
8. The tower $L_0 \subseteq L \subseteq E$\\
9. Deformation of $V$-modules\\
10. Weak Jacobi forms \\
11. The quadratic space $(E, \langle \ , \ \rangle)$\\
12. Commutants\\
13. $U$-weights\\
14. Lattice subalgebras of $V$\\
15. The tripartite subVOA of $V$\\
16. The invariants $\tilde{c}$ and $l$\\

\bigskip

\noindent
1. {\bf Statement of main results}

\medskip
This paper concerns the algebraic structure of  \emph{strongly regular} vertex operator algebras (VOAs). A VOA $V= (V, Y, \mathbf{1}, \omega)$ is called \emph{regular} if it is \emph{rational} (admissible $V$-modules are semisimple) and \emph{$C_2$-cofinite} (the span of $u(n)v\ (u, v \in V, n \leq -2)$ has finite codimension in $V$). It is strongly regular if, in addition, the $L(0)$-grading (or \emph{conformal grading}) given by
$L(0)$-weight has the form
\begin{eqnarray}\label{specdecomp1}
V = \mathbb{C}\mathbf{1} \oplus V_1 \oplus \hdots
\end{eqnarray}
and all states in $V_1$ are \emph{quasiprimary} (i.e.\ annihilated by $L(1)$).
Apart from the still-undecided question of the relationship between  rationality and $C_2$-cofiniteness, changing any of the assumptions in 
the definition of strong regularity will result in VOAs with quite different properties (cf.\ \cite{DM5}). 
Such VOAs are of interest in their own right, but we will not deal with them here.

\medskip
To describe the main results, we need some basic facts about strongly regular VOAs $V$ that will be assumed here and reviewed in more detail in later Sections. $V$ is equipped with an essentially unique 
nonzero, invariant, bilinear form $\langle \ , \ \rangle$, and $V$ is simple if, and only if, $\langle \ , \ \rangle$ is \emph{nondegenerate}. We assume this is the case from now on. Then $V_1$ carries the structure of a \emph{reductive} Lie algebra and all Cartan subalgebras of $V_1$ (maximal (abelian) toral Lie subalgebras)
are conjugate in Aut$(V)$. We also refer  them as \emph{Cartan subalgebras of $V$}. We say that a subspace $U \subseteq V$ is \emph{nondegenerate} if the restriction of
$\langle \ , \ \rangle$ to $U \times U$ is nondegenerate. For example, the Cartan subalgebras of $V$
and the solvable radical of $V_1$ are nondegenerate. We refer to the dimension of $H$ as the \emph{Lie rank} of $V$.

\medskip
A \emph{subVOA} of $V$ is a subalgebra $W = (W, Y, \mathbf{1}, \omega')$ with a conformal vector $\omega'$ that may not coincide with the conformal vector $\omega$ of $V$. If $\omega = \omega'$ we say that $W$ is a \emph{conformal} subVOA. $V$ contains a unique minimal conformal subVOA (with respect to inclusion), namely the \emph{Virasoro subalgebra} generated by $\omega$.
A basic example of a subVOA is the Heisenberg theory
$(M_U, Y, \mathbf{1}, \omega_U)$ generated by a nondegenerate subspace $U$ of a Cartan subalgebra of $V$.
$M_U$ has rank (or central charge) $\dim U$ and conformal vector
\begin{eqnarray}\label{moreomegadefs}
\omega_U:= 1/2\sum_i h^i(-1)h^i,
\end{eqnarray}
(for any orthonormal basis $\{h^i\}$ of $U$).
A \emph{lattice theory} is a VOA
$V_L$ corresponding to a positive-definite, even lattice $L$. 

\medskip
We can now state the main result. 

\medskip
\noindent
{\bf Theorem 1}: Let $V$ be a strongly regular, simple VOA, and suppose that $U \subseteq H \subseteq V$ where $H$ is a Cartan subalgebra of $V$ and $U$ is a nondegenerate subspace.
Let $\omega_U$ be as in (\ref{moreomegadefs}). Then the following hold: 
\begin{eqnarray}
&&(a)\ \mbox{There is a unique maximal subVOA  
$W \subseteq V$ with} \notag\\
&&\ \ \ \ \  \mbox{conformal vector $\omega_U$.} \notag\\
&&(b)\ \mbox{$W \cong V_{\Lambda}$ is a lattice theory, where $\Lambda \subseteq U$ is a} \notag \\
&&\ \ \ \ \ \mbox{positive-definite even lattice with  $\dim U = \mbox{rk} \Lambda$}.  \label{Widentify}
\end{eqnarray}

\medskip \noindent
{\bf Remark 2}. 
1. Part (a) - the fact that there is a unique maximal subVOA $W$ with conformal vector
$\omega_U$ - is elementary; it follows from the theory of commutants \cite{FZ} (cf.\ Section 12). 
The main point of the Theorem is  the \emph{identification} of $W$ as a lattice theory.\\
 \noindent
2. $U$ is a Cartan subalgebra of $W$. Thus, every nondegenerate subspace of $H$ is a Cartan
subalgebra of a lattice subVOA of $V$.

\bigskip
Theorem 1 has many consequences. We discuss some of them here, 
deferring a fuller discussion until later Sections.
 We can apply Theorem 1 with $U=$ rad$(V_1)$, and this leads to the next result.
  
\medskip \noindent
{\bf Theorem 3}: Suppose that $V$ is a strongly regular, simple $VOA$. There is a \emph{canonical} conformal subVOA
\begin{eqnarray}\label{tripdef}
T = W \otimes G \otimes Z,
\end{eqnarray}
the tensor product of subVOAs $W, G, Z$ of $V$ with the following properties: \\
(a)\ $W \cong V_{\Lambda}$ is a lattice theory and
$\Lambda$ has minimal length at least $4$; \\
(b)\ $G$ is the tensor product of affine Kac-Moody algebras of positive integral level; \\
(c)\ $Z$ has no nonzero states of weight $1$: 
$Z = \mathbb{C}\mathbf{1} \oplus Z_2 \oplus \hdots$

\medskip \noindent
{\bf Remark 4}.
 The gradings on $W, G$ and $Z$ are compatible with that on $V$ in the sense that the $n^{th}$ graded piece of each of them is contained in $V_n$. $T$ has the tensor product grading, and in particular $T_1 = W_1 \oplus G_1 = V_1$. Indeed, $W_1=$ rad$(V_1)$ and $G_1$ is the Levi factor of $V_1$.  Thus the weight $1$ piece of
$V$ is contained in a rational subVOA of standard type, namely a tensor product of a lattice theory and affine Kac-Moody algebras.

\medskip
To a certain extent, Theorem 3 reduces the study of strongly regular VOAs to the following:
(A) proof that $Z$ is strongly regular; (B) study of strongly regular VOAs with no nonzero weight $1$
states; (C) extension problem for strongly regular VOAs, i.e.\ characterization of
 the strongly regular VOAs that contain a \emph{given} strongly regular conformal subVOA $T$.
 For example, we have the following immediate consequence of Theorem 3 and Remark 4. 
 
 \medskip \noindent
{\bf Theorem 5}: Suppose that $V$ is a strongly regular, simple $VOA$ such
that the conformal vector $\omega$ lies in the subVOA $\langle V_1 \rangle$ generated
by $V_1$.  Then the canonical conformal subalgebra (\ref{tripdef}) is a rational subVOA
\begin{eqnarray*}
T = W \otimes G,
\end{eqnarray*}
where $W$ and $G$ are as in the statement of Theorem 2.

\medskip \noindent
{\bf Remark 6}. 
 Let $\mathcal{C}$ consist of the (isomorphism classes of) VOAs satisfying the assumptions
 of the Theorem. $\mathcal{C}$ contains all lattice theories, all simple affine Kac-Moody 
 VOAs of positive integral level (Siegel-Sugawara construction), and it is 
 closed
 with respect to tensor products and extensions in the sense of (C) above.
  Theorem 5 says that
every VOA in $\mathcal{C}$ arises this way, i.e.\ an extension of a tensor product of a lattice theory and
affine Kac-Moody theories.

 \medskip
There are applications of Theorem 1 
 to inequalities involving the Lie rank $l$ and the \emph{effective central charge} $\tilde{c}$ of $V$. These lead to characterizations of some classes of strongly rational VOAs $V$ according to these invariants. For example, we have

 \medskip \noindent
{\bf Theorem 7}: Let $V$ be a strongly regular, simple  VOA of effective central charge $\tilde{c}$ and Lie rank $l$. The following are equivalent:
\begin{eqnarray*}
&&(a)\ \tilde{c} < l+1, \\
&&(b)\ V\ \mbox{contains a conformal subalgebra isomorphic to a} \\
&&\ \ \ \ \ \mbox{tensor product $V_{\Lambda} \otimes L(c_{p, q}, 0)$ of a lattice theory of rank $l$}\\
&&\ \ \ \  \ \mbox{and a \emph{simple Virasoro VOA} in the \emph{discrete series}.}
\end{eqnarray*}

\medskip \noindent
{\bf Remark 8.}  1. We \emph{always} have
$l \leq \tilde{c}$ (\cite{DM1}). \\
2. Define a
\emph{minimal model} as a strongly regular simple VOA whose Virasoro subalgebra lies in the 
\emph{discrete series}. 
 The case $l=0$ of Theorem 7 characterizes minimal models as those
 strongly regular simple VOAs which have $\tilde{c}<1$. This is, of course,
 very similar to
 the classification of minimal models  in physics (cf.\ \cite{FMS}, Chapters 7 and 8), where attention
 is usually restricted to the \emph{unitary} case, where $c=\tilde{c}$, or equivalently $q=p+1$,  in the notation of Theorem 7.
 Minimal models with $\tilde{c}=c$ were treated rigorously in \cite{DW2};
 our approach allows us to remove any assumptions about $c$ and permits $l$ to be nonzero.

\medskip
Our results  continue, and in some cases complete, lines of thought in \cite{DM1} and \cite{DM2} having to do with the weight $1$ subspace $V_1$ of $V$ and its embedding in $V$. These include the invariant bilinear form of $V$, the nature of the Lie algebra of $V_1$ and its action on $V$-modules, automorphisms of $V$ induced by exponentiating weight one states, deformations of $V$-modules using weight one states, and (more recently \cite{KM}) 
 weak Jacobi form trace functions defined by weight one states. A fuller account might also have included simple currents arising from deformations by weight one states \cite{DLM4}, although we do not treat this subject here.\
 These topics constitute a very satisfying Chapter in the theory of rational VOAs, and the Heidelberg Conference presented itself as a great opportunity to review this set of ideas. I am grateful to the organizers, Professors Winfried Kohnen and Rainer Weissauer, for giving me the chance to do so.

\bigskip
\noindent
2. {\bf Background}

\medskip
A vertex operator algebra (VOA) is a quadruple $(V, Y, \mathbf{1}, \omega)$, often denoted simply by
$V$, satisfying the usual axioms. For these and other background results in VOA theory, we refer the reader to \cite{LL}. We write vertex operators as
\begin{eqnarray*}
Y(v, z) &=& \sum_{n \in \mathbb{Z}} v(n)z^{-n-1} \ \ (v \in V), \\
Y(\omega, z) &=& \sum_{n \in \mathbb{Z}} L(n)z^{-n-2}.
\end{eqnarray*}
Useful identities that hold for all $u, v \in V, p, q \in \mathbb{Z}$ include
\begin{eqnarray}
[u(p), v(q)] &=& \sum_{i=0}^{\infty} {p \choose i}(u(i)v)(p+q-i), \label{commform}\\
\{u(p)v\}(q) &=& \sum_{i=0}^{\infty} (-1)^i{p \choose i}(u(p-i)v(q+i) - (-1)^pv(q+p-i)u(i)), \notag
\end{eqnarray}
called the \emph{commutator formula} and \emph{associativity  formula} respectively.

\medskip
We assume throughout that $V$ is a \emph{simple} VOA that is strongly regular as defined in Section 1. 
One of the main consequences of rationality is the fact that, up to isomorphism, 
there are only \emph{finitely many} ordinary irreducible $V$-modules (\cite{DLM2}). We let
$\mathcal{M} :=\{(M^1, Y^1), \hdots, (M^r, Y^r)\}$ denote this set, with $(M^1, Y^1)=(V, Y)$. It is conventional to use $u(n)$ to denote the $n^{th}$ mode of $u \in V$ acting on any $V$-module, the meaning usually being clear from the context, however it will sometimes be convenient to distinguish
some of these modes. In particular, we often write $Y^j(u, z) := \sum_{n \in \mathbb{Z}}u_j(n)z^{-n-1} \ (u \in V), Y^j(\omega, z):= \sum_{n \in \mathbb{Z}} L_j(n)z^{-n-1}$, dropping the index $j$ from the notation when $j=1$.

\bigskip
\noindent
3. {\bf Invariant bilinear form}

\medskip
An \emph{invariant bilinear form} on $V$ is a bilinear map
$\langle \ , \ \rangle: V \times V \rightarrow \mathbb{C}$ satisfying
\begin{eqnarray}\label{invformdef1}
\langle Y(a, z)b, c \rangle = \langle b, Y(e^{zL(1)}(-z^{-2})^{L(0)}a, z^{-1} \rangle\ \ (a, b, c \in V).
\end{eqnarray}
Such a form is necessarily \emph{symmetric} ([FHL], Proposition 5.3.6).

\medskip
A theorem of Li \cite{L1} says that  there is a linear isomorphism 
between $V_0/L(1)V_1$ and the space of  invariant bilinear forms on $V$.
Because $V$ is strongly regular then $V_0/L(1)V_1 = \mathbb{C}\mathbf{1}$, so a nonzero
invariant bilinear form exists and it is \emph{uniquely determined} up to scalars.

\medskip
 If $a \in V_k$ is quasi-primary
 then (\ref{invformdef1}) says that
\begin{eqnarray}\label{invformdef2}
\langle a(n)b, c \rangle = (-1)^k\langle b, a(2k-n-2)c \rangle \ \ \ (n \in \mathbb{Z}).
\end{eqnarray}
In particular, this applies if $a \in V_1$ (because $V$ is assumed to be strongly regular), or if 
$a = \omega$ is the conformal vector ($\omega$ is \emph{always} quasiprimary). First apply
(\ref{invformdef2}) with $a = \omega$ and $n=1$, noting that
$\omega(1)=L(0)$. Then $k=2$ and we obtain
\begin{eqnarray*}
\langle L(0)b, c \rangle = \langle b, L(0)c \rangle.
\end{eqnarray*}
It follows that eigenvectors of $L(0)$ with distinct eigenvalues are necessarily perpendicular with respect to $\langle\ , \ \rangle$.
Thus (\ref{specdecomp1}) is an orthogonal direct sum
 \begin{eqnarray}\label{specdecomp2}
V &=& \mathbb{C}\mathbf{1} \perp V_1 \perp \hdots  
\end{eqnarray}
(Here and below, for subsets $A, B \subseteq V$ we write $A \perp B$ if $\langle a, b \rangle = 0$ for all
 $a \in A, b \in B$.)
 
 \medskip
 The radical $R$ of $\langle \ , \ \rangle$ is an \emph{ideal}. Because we are assuming
 that $V$ is simple then $R=0$ and $\langle \ , \ \rangle$ is \emph{nondegenerate}.
 In particular, $\langle \mathbf{1}, \mathbf{1} \rangle \not= 0$.
 In what follows, we fix the form so that 
\begin{eqnarray}\label{normchoice}
\langle \mathbf{1}, \mathbf{1} \rangle = -1.
\end{eqnarray}
Note also that by (\ref{specdecomp2}),
 the restriction of $\langle\ , \ \rangle$ to each $V_n \times V_n$ is also nondegenerate.

\bigskip
\noindent
4. {\bf The Lie algebra on $V_1$}

\medskip
The bilinear product $[uv] := u(0)v\ (u, v \in V_1)$ equips $V_1$ with the structure of a Lie algebra. Applying (\ref{invformdef2}) with  $u, v \in V_1$, we obtain
$\langle u, v \rangle = \langle u(-1)\mathbf{1}, v \rangle = -\langle \mathbf{1}, u(1)v \rangle$. With the convention (\ref{normchoice}), it follows that
\begin{eqnarray*}
u(1)v = \langle u, v \rangle \mathbf{1} \ \ (u, v \in V_1).
\end{eqnarray*}

Because $V$ is strongly regular, a theorem of Dong-Mason \cite{DM1} says that the Lie algebra 
on $V_1$ is \emph{reductive}. (This result is discussed further in Section 7 below.)
So there is a canonical decomposition
$V_1 = A \perp S$
where $A=$ Rad$(V_1)$ is an \emph{abelian} ideal and
$S$ is the (semisimple) Levi factor. The decomposition of $S$ into a direct sum of simple Lie algebras
$\oplus_i \frak{g}_i$ is also an orthogonal sum with respect to $\langle\ , \ \rangle$. 

\medskip
There is a refinement of this decomposition, established in \cite{DM2}, namely
\begin{eqnarray}\label{LAdecomp}
V_1 = A \perp \frak{g}_{1, k_1} \perp \hdots \perp \frak{g}_{s, k_s}
\end{eqnarray}
where $\frak{g}_{i}$ is a simple Lie algebra and $k_i$ is a positive integer (the \emph{level}).  

\medskip
To explain what this means, for $U \subseteq V$ let $\langle U \rangle$ be the 
\emph{subalgebra} of $V$ 
\emph{generated} by $U$. $\langle U \rangle$ is spanned by states
$u = u_1(n_1)\hdots u_t(n_t)\mathbf{1}$ with $u_1, \hdots, u_t \in U,\ n_1, \hdots, n_t \in \mathbb{Z}$,
and equipped with vertex operators  defined as the \emph{restriction}
of $Y(u, z)$ to $\langle U \rangle$. 

\medskip
 It is proved in \cite{DM2} that
there is an isomorphism of VOAs
$\langle \frak{g}_{i} \rangle \cong L_{\frak{g}_i}(k_i, 0)$, where $L_{\frak{g}_i}(k_i, 0)$ is the 
simple VOA (or WZW model) corresponding to the affine Lie algebra 
$\widehat{\frak{g}_i}$ determined by 
$\frak{g}_i$, of positive integral level $k_i$.
Orthogonal Lie algebras in (\ref{LAdecomp}) determine mutually commuting WZW models. So the meaning of (\ref{LAdecomp}) is that the canonical subalgebra $G$ of
$V$ generated by $S$  satisfies
\begin{eqnarray}\label{Wiso}
G \cong L_{\frak{g}_1}(k_1, 0) \otimes \hdots \otimes L_{\frak{g}_s}(k_s, 0).
\end{eqnarray}
In particular, $G$ is a rational VOA  equipped with the canonical conformal vector 
$\omega_G$ arising from the Sugawara construction associated to each tensor factor
(\cite{FZ}, \cite{LL}).

\medskip
Because $V_1$ is reductive, it has a \emph{Cartan subalgebra}, that is a maximal (abelian) toral subalgebra, and all Cartan subalgebras are conjugate in Aut$(V_1)$. (See the following Section for further discussion.)
Let $H \subseteq V_1$ be a Cartan subalgebra of $V_1$, say of rank $l$. By Lie theory, the restriction of $\langle \ , \ \rangle$ to $H \times H$ is nondegenerate. We also call $H$ a \emph{Cartan subalgebra of $V$}.

\bigskip \noindent
5. {\bf Automorphisms}

\medskip
An \emph{automorphism} of $V$ is an invertible linear map
$g: V \rightarrow V$ such that $g(\omega) = \omega$ and $ga(n)g^{-1}=g(a)(n)$ for all $a \in V, n \in \mathbb{Z}$, i.e.
\begin{eqnarray}\label{autact}
gY(a, z)g^{-1} = Y(g(a), z).
\end{eqnarray}
The set of all automorphisms is a group Aut$(V)$.  Because $g\omega(n)g^{-1} = g(\omega)(n)
=\omega(n)$, it follows in particular that $g$ commutes with $L(0) = \omega(1)$. Therefore, Aut$(V)$ acts on each $V_n$. The uniqueness of $\langle \ , \ \rangle$ implies that 
\begin{eqnarray}\label{Gforminv}
\mbox{Aut$(V)$ leaves $\langle \ , \ \rangle$ invariant}.
\end{eqnarray}
So each $V_n$ affords an \emph{orthogonal} representation of Aut$(V)$.

\medskip
One checks (e.g.\ using induction and  (\ref{commform})) that for $n \geq 0$, 
\begin{eqnarray*}
(u(0)^nv)(q) = \sum_{i=0}^n (-1)^i{ n \choose i}u(0)^{n-i}v(q)u(0)^{i} \ \ (u, v \in V, \ q \in \mathbb{Z}).
\end{eqnarray*}
Therefore,
\begin{eqnarray*}
\left(e^{u(0)}v\right)(q) &=&\sum_{n=0}^{\infty} \frac{1}{n!} (u(0)^nv)(q) \\
&=&\sum_{n = 0}^{\infty}\sum_{i=0}^n \frac{(-1)^i}{i!(n-i)!} u(0)^{n-i}v(q)u(0)^{i} \\
&=&e^{u(0)}v(q)e^{-u(0)},
\end{eqnarray*}
showing that (\ref{autact}) holds with $g = e^{u(0)}$. If we further assume that $u \in V_1$
then we obtain using (\ref{commform}) that
\begin{eqnarray*}
u(0)\omega 
&=& -[\omega(-1), u(0)]\mathbf{1} = - \sum_{i=0}^{\infty} (-1)^i (\omega(i)u)(-1-i)\mathbf{1}  \\
&=& - ((L(-1)u)(-1) -(L(0)u)(-2))\mathbf{1} = 0.
\end{eqnarray*}

It follows  that $\{e^{u(0)} \ | \ u \in V_1\}$ is a set of automorphisms of $V$. Let
\begin{eqnarray*}
\frak{G} = \langle e^{u(0)} \ | \ u \in V_1 \rangle
\end{eqnarray*}
be the group they generate.
It is clear from the classical relation between Lie groups and Lie algebras that $\frak{G}$ is the adjoint form of the complex Lie group  associated with $V_1$. So there is a containment
\begin{eqnarray*}
\frak{G} \unlhd \mbox{Aut}(V).
\end{eqnarray*}
(Normality holds because
if $g \in$ Aut$(V)$ and $u \in V_1$ then $g(u)\in V_1$ and $ge^{u(0)}g^{-1}= e^{gu(0)g^{-1}}
=e^{g(u)(0)}$.)

\medskip
One consequence of this is the following. Because $\frak{G}$  acts transitively on the set of Cartan subalgebras of $V_1$, it follows \textit{ipso facto} that Aut$(V)$ also acts transitively on the set of Cartan subalgebras of $V_1$ (or of $V$). Thus the choice of a Cartan subalgebra in $V$ is unique up to automorphisms of $V$, in parallel with the usual theory of semisimple Lie algebras.

\bigskip
\noindent 
6. {\bf Projective action of Aut$(V)$ on $V$-modules}

\medskip
There is a natural action of Aut$(V)$ on the set $\mathcal{M}$ of (isomorphism classes of) irreducible $V$-modules
$\{(M^j, Y^j) \ | \ 1 \leq j \leq r\}$ \cite{DM3}. Briefly, the argument is as follows.
For $g \in$ Aut$(V)$ and an index $j$, one checks that the pair $(M^j, Y^j_g)$ defined by
$Y^j_g(v, z):= Y^j(gv, z)\ (v \in V)$ is itself an irreducible $V$-module. Thus the action of Aut$(V)$ on
$\mathcal{M}$ is defined by $g: (M^j, Y^j) \mapsto (M^j, Y^j_g)$.

\medskip
Because $\mathcal{M}$ is finite and $\frak{G}$ is connected, the action of $\frak{G}$ is necessarily 
\emph{trivial}. Hence, if we fix the index $j$, then for $g \in$ Aut$(V)$ there is an isomorphism of
$V$-modules $\alpha_g: (M^j, Y^j) \mapsto (M^j, Y^j_g)$, i.e.
\begin{eqnarray}\label{gact1}
\alpha_g Y^j(u, z) = Y^j_g(u, z)\alpha_g = Y^j(gu, z)\alpha_g\ \ (u \in V).
\end{eqnarray}
Because $M^j$ is irreducible, $\alpha_g$ is uniquely determined up to an overall nonzero scalar
(Schur's Lemma).

\medskip
When $j=1$, so that $M^j=V, \alpha_g$ coincides with $g$ itself, the scalar being implicitly
determined
by the additional condition $g(\omega)=\omega$. Generally, we find from (\ref{gact1}) that
\begin{eqnarray*}
\alpha_{gh}Y^j(u, z)\alpha_{gh}^{-1} = \alpha_g \alpha_hY^j(u, z)\alpha_h^{-1}\alpha_g^{-1}
\ (g, h \in \mbox{Aut}(V)),
\end{eqnarray*}
so that by Schur's Lemma once more there are scalars $c_j(g, h)$ satisfying
\begin{eqnarray*}
\alpha_{gh} = c_j(g, h)\alpha_g \alpha_h.
\end{eqnarray*}
The map $c_j: \frak{G}\times\frak{G} \rightarrow \mathbb{C}^*, \ (g, h) \mapsto c_j(g, h),$
is a $2$-cocycle on $\frak{G}$. It defines a projective action $g \mapsto \alpha_g$
of $\frak{G}$ on $M^j$ that satisfies (\ref{gact1}).

\medskip
While the projective action of $\frak{G}$ on $M^1=V$ reduces to the linear action previously considered, the $2$-cocycles $c_j$ are generally \emph{nontrivial}, i.e.\ they are not $2$-coboundaries. A well-known example is the VOA $V:= L_{\frak{sl}_2}(1, 0)$, i.e.\
the level $1$ WZW model of type $sl_2$, which is isomorphic to the lattice theory
$V_{A_1}$ defined by the $A_1$ root lattice. In this case we have
$V_1= \frak{sl}_2$, with adjoint group $\frak{G}= SO_3(\mathbb{R})$. There are just two irreducible
$V$-modules, corresponding to the two cosets of $A_1$ in its dual lattice $A_1^*:=(1/\sqrt{2})A_1$,
and their direct sum is the Fock space for the super VOA $V_{A_1^*}$. The automorphism group of this
SVOA is $SU_2(\mathbb{C})$, in other words the projective action of
$\frak{G}$ on $M^2$ lifts to a linear action of its proper $2$-fold (universal) covering group.

\bigskip
\noindent
7. {\bf Complete reducibility of the $V_1$-action}

\medskip
We discuss the following result.
\begin{eqnarray}\label{crH}
&&\mbox{The Lie algebra $V_1$ is \emph{reductive}, and its action on each} \notag\\
&&\mbox{simple $V$-module $(M^j, Y^j)$ is \emph{completely reducible}.}
\end{eqnarray}
This follows from results in \cite{DM1} and \cite{DG}. We will need some of the details later, so we sketch the proof.

\medskip
Each irreducible $V$-module $M^j$ has a direct sum decomposition into finite-dimensional
$L_j(0)$-eigenspaces
\begin{eqnarray}\label{M^jgrading}
M^j = \bigoplus_{n= 0}^{\infty} M^j_{n+\lambda_j},
\end{eqnarray}
where $\lambda_j$ is a constant called the \emph{conformal weight} of $M^j$.
Each $M^j_{n+\lambda_j}$ is a module for the Lie algebra $V_1$, acting by the zero mode
$u_j(0) \ (u \in V_1)$, and (\ref{crH}) amounts to the assertion that each of these actions is completely reducible. The simple summands $\frak{g}_i\ (1 \leq i \leq s)$ of $V_1$ act completely reducibly by 
Weyl's theorem, so the main issue is to show that the abelian radical $A$ of $V_1$ (cf.\ (\ref{LAdecomp})) acts semisimply.

\medskip
The first step uses a formula of Zhu
\cite{Z}. The case we need may be stated as follows (cf.\ \cite{DM1}):
\begin{eqnarray}\label{Zmodinv1}
&&\ \ \ \ \ \ \ \ \mbox{Suppose that $u, v \in V_1$. Then for $1 \leq j \leq r,$}  \\
&& \mbox{Tr}_{M^j} u_j(0)v_j(0)q^{L(0)-c/24} = Z_{M^j}(u[-1]v, \tau) - \langle u, v \rangle E_2(\tau)Z_{M^j}(\tau). \notag
\end{eqnarray}
The notation, which is standard, is as follows (\cite{Z}, \cite{DLM3}): for $w \in V$, \\
$Z_{M^j}(w, \tau):=$ 
Tr$_{M^j}o_j(w)q^{L(0)-c/24}$ is the graded trace of the \emph{zero} mode $o_j(w)$ for the action of 
$w$ on $M^j, u[-1]$ is the $-1^{st}$ \emph{square bracket} mode for $u$,
and $E_2(\tau) = -1/12 +2\sum_{n=1}^{\infty} \sum_{d|n} dq^n$ is the usual weight 
$2$ Eisenstein series.

\medskip
Next we show that if $\langle u, v \rangle \not= 0$ then for some index $j$ we have
\begin{eqnarray*}
Z_{M^j}(u[-1]v, \tau) \not= \langle u, v \rangle E_2(\tau)Z_{M^j}(\tau).
\end{eqnarray*}
 Indeed, if this does not hold, we can obtain a contradiction  using Zhu's modular-invariance theorem \cite{Z} and the exceptional transformation law for $E_2(\tau)$ (cf.\ \cite{DM1}, Section 4 for details). From (\ref{Zmodinv1})  we can conclude that if $\langle u, v \rangle \not= 0$ then there is an index $j$
 such that
\begin{eqnarray}\label{Zmodinv2}
 \mbox{Tr}_{M^j} u_j(0)v_j(0) \not= 0.
 \end{eqnarray}

Now suppose that $u \in V_1$ lies in the nil radical of $V_1$. Then $u(0)$ annihilates every
simple $V_1$-module, and in particular the lhs of (\ref{Zmodinv2}) necessarily vanishes for each $j$.
Therefore $\langle u, v \rangle = 0\ (v \in V_1)$, whence $u=0$. This shows that $V_1$
is indeed reductive. 

\medskip
It is known that a  VOA is \emph{finitely generated} (f.g.)\ if it is $C_2$-cofinite \cite{GN}, \cite{B},  or if it is
rational \cite{DW1}. So certainly a regular VOA is f.g.
We need this mainly because Griess and Dong proved \cite{DG} that the automorphism group of a f.g.\ VOA is a (complex) \emph{algebraic group}. It follows that the subgroup
$\frak{A} \unlhd \frak{G}$ generated by the exponentials $e^{u(0)}\ (u \in A=$ rad($V_1)$) is an abelian  algebraic subgroup, and that $A$ itself is the direct sum of two Lie subalgebras corresponding to the
unipotent and semisimple parts of $\frak{A}$. By the same argument as above, the Lie subalgebra corresponding to the unipotent part necessarily vanishes, so that $\frak{A}$ is a complex torus
and $A$ consists of semisimple operators. In particular, (\ref{crH}) holds.

\medskip
(\ref{crH}) was first stated in \cite{DM1}, although the proof there is incomplete. It would be interesting
to find a proof that does not depend on  the theory of algebraic groups.

\bigskip
\noindent
8. {\bf The tower $L_0 \subseteq L \subseteq E$}

\medskip
Fix a Cartan subalgebra $H \subseteq V_1$ of rank $l$, say. We have seen in Section 7
that all of the operators
$u_j(0)\ (u \in H, 1 \leq j \leq r)$ are semisimple. We set 
\begin{eqnarray}\label{Edef}
E &=& \{ u \in H \ | \ u(0) \ \mbox{has eigenvalues in}\  \mathbb{Q}\},  \notag\\
L &=& \{ u \in H \ | \ u_j(0) \ \mbox{has eigenvalues in}\ \mathbb{Z},\ 1 \leq j \leq r \},\\
L_0 &=& \{ u \in H \ | \ u(0) \ \mbox{has eigenvalues in}\ \mathbb{Z} \}. \notag
\end{eqnarray}
$E$ is a $\mathbb{Q}$-vector space in $H$ and $L \subseteq L_0 \subseteq E$ are additive subgroups.  

\medskip
Let $\frak{H} \subseteq \frak{G}$ be the group generated by exponentials
 $e^{2\pi i u(0)} \ (u \in H)$. From Section 7, $\frak{H}$ is a complex torus
 $\frak{H} \cong (\mathbb{C}^*)^l$. There is a short exact sequence
\begin{eqnarray*}
0 \rightarrow L_0/L \rightarrow H/L \stackrel{\varphi}{\rightarrow} \frak{H} \rightarrow 1
\end{eqnarray*}
 where $\varphi$ arises from the morphism $u \mapsto e^{2\pi i u(0)}\ (u \in H)$.
 $H/L$ is the covering group of $H/L_0$ that acts linearly on each irreducible module $M^j$ as described in Section 6,
 and $H/L_0 \cong \frak{H}$.
 
 \medskip
Because $V$ is f.g. there is an integer $n_0$ such that $V = \langle \oplus_{n=0}^{n_0}V_n\rangle$. Then $e^{2\pi i u(0)}\ (u\in H)$ is the identity if, and only if, its restriction to $\oplus_{n=0}^{n_0} V_n$ is the identity. It follows  that the eigenvalues of
$u(0)$ for $u \in E$ have bounded denominator, whence 
\begin{eqnarray}\label{Etor}
E/L_0 = \mbox{Torsion}(H/L_0) \cong (\mathbb{Q}/\mathbb{Z})^l. 
\end{eqnarray}
In particular, $E$ contains a $\mathbb{C}$-basis of $H$.

\bigskip
\noindent
 9. {\bf Deformation of $V$-modules}

\medskip
In \cite{L2}, Proposition 5.4, Li showed how to deform (twisted) $V$-modules using a certain operator 
$\Delta(z)$. We describe the special case that we need here. 
See \cite{KM} for further details of the calculations below, and \cite{DLM4} for 
further development of the theory.

\medskip
Fix $u \in L_0$ (\ref{Edef}), and set
\begin{eqnarray*}
\Delta_{u}(z) &:=& z^{u(0)}\exp\left\{- \sum_{k \geq 1} \frac{u(k)}{k}(-z)^{-k} \right\}.
\end{eqnarray*}
For an irreducible $V$-module $(M^{j'}, Y^{j'})$, set
\begin{eqnarray*}
Y^{j'}_{\Delta_{u}(z)}(v, z) &:=& Y^{j'}(\Delta_{u}(z)v, z) \ \ (v \in V).
\end{eqnarray*}
Because $u(0)$ has eigenvalues in $\mathbb{Z}$ then $e^{2\pi i u(0)}$ is the identity automorphism of $V$. In this case, Li's result says that  there is an isomorphism of $V$-modules
\begin{eqnarray}\label{Liiso}
 (M^{j'}, Y^{j'}_{\Delta_{u}(z)}) \cong (M^j, Y^j)
\end{eqnarray}
for some $j$. (Technically, Li's results deal with \emph{weak} $V$-modules. In the case that we are dealing with, when $V$ is regular, the results apply to ordinary irreducible $V$-modules, as stated.)
Thus there is a linear isomorphism $\psi: M^{j'} \stackrel{\cong}{\rightarrow} M^j$ satisfying
\begin{eqnarray}\label{psiintertwine}
\psi^{-1}Y^j(v, z)\psi = Y^{j'}(\Delta_u(z)v, z) \ \ (v \in V).
\end{eqnarray}

\medskip
 In (\ref{psiintertwine}) we choose $j'=1$ (so $(M^{j'}, Y^{j'}) = (V, Y)$), $v = \omega$, 
 and apply both sides to $\mathbf{1}$. We obtain after some calculation that
 \begin{eqnarray}\label{calc1}
\psi^{-1} L^j(0)\psi(\mathbf{1}) = 1/2\langle u, u \rangle \mathbf{1}.
\end{eqnarray}
The $L(0)$-grading on $M^j$ is described in (\ref{M^jgrading}). If $\psi(\mathbf{1}) = \sum_n a_n$
with $a_n \in M^j_{n+\lambda_j}$ then
$1/2\langle u, u \rangle\sum_n a_n = \sum_n (n+\lambda_j)a_n$. This shows that
$\psi(\mathbf{1}) \in M^j_{n_0+\lambda_j}$ for some integer $n_0$, and moreover
\begin{eqnarray}\label{uurat1}
1/2\langle u, u \rangle = n_0+\lambda_j.
\end{eqnarray}

We use (\ref{uurat1}) in conjunction with another
Theorem (\cite{DLM3}, \cite{AM}) that says that (for regular $V$) the conformal weight
$\lambda_j$ of the irreducible $V$-module $M^j$ lies in $\mathbb{Q}$. Then it is immediate from  (\ref{uurat1})  that $\langle u, u \rangle \in \mathbb{Q}.$ The only condition on $u$ here is that $u \in L_0$. Because $E/L_0$ is a torsion group
(\ref{Etor}) we obtain
\begin{eqnarray}\label{lambdacong}
\langle u, u \rangle \in \mathbb{Q} \ \ (u \in E).
\end{eqnarray}

\medskip
Arguing along similar lines,  we can also prove the following: (i) if $u \in E$, all eigenvalues of
the operators $u_j(0)$ lie in $\mathbb{Q}\ (1 \leq j \leq r)$; (ii) if $u \in L_0$ then the denominators of the eigenvalues
of $u_j(0)$ divide the l.c.m.\ $M$ of the denominators of the conformal weights $\lambda_j$.
In other words, $L_0/L$ is a torsion abelian group of exponent dividing $M$. (It is also f.g., as we shall see. So $L_0/L$  is actually a finite abelian group.)

\bigskip
\noindent
10. {\bf Weak Jacobi forms}

\medskip
The paper \cite{KM} develops an extension of Zhu's theory of partition functions \cite{Z} to the context of 
 \emph{weak Jacobi forms}. We discuss background sufficient for our purposes. For the general theory of Jacobi forms, cf.\ \cite{EZ}.

\medskip
We continue with a strongly regular VOA $V$. Let $h \in L$ (cf.\ (\ref{Edef})). For 
$j$ in the range $1 \leq j \leq r$, define
\begin{eqnarray*}
J_{j, h}(\tau, z) := \mbox{Tr}_{M^j} q^{L_j(0)-c/24}\zeta^{h_j(0)},
\end{eqnarray*}
where $c$ is the central charge of $V$. (The definition makes sense because we 
have seen that $h_j(0)$ is a semisimple operator.) Notation is as follows:
$q:=e^{2\pi i \tau}, \zeta:=e^{2\pi i z}, \tau \in \mathbb{H}$ (complex upper half-plane), $z \in \mathbb{C}$.
The main result \cite{KM} is that
$J_{j, h}(\tau, z)$ is holomorphic in $\mathbb{H}\times\mathbb{C}$ and satisfies
the following functional equations for 
all $\gamma=\left(\begin{array}{cc}a&b \\c&d\end{array}\right) \in SL_2(\mathbb{Z}),
(u, v)\in \mathbb{Z}^{2}, 1 \leq i \leq r$:
\begin{eqnarray}
&&\ (i)\ \mbox{there are scalars $a_{ij}(\gamma)$ depending only on $\gamma$ such that} \notag\\
&&\ \ \ \ J_{i, h} \left(\gamma\tau, \frac{z}{c\tau+d}\right) =
e^{\pi i cz^2\langle h, h \rangle/(c\tau+d)}\sum_{j=1}^r a_{ij}(\gamma) \label{Gammaact}
J_{j, h}(\tau, z), \\
&&(ii)\ \mbox{there is a permutation $j \mapsto j'$ of $\{1, \hdots, r\}$ such that} \notag\\
&&\ \ \ \ \ J_{j, h}(\tau, z+u\tau+v) 
=e^{-\pi i \langle h, h \rangle(u^2 \tau +2uz)} \label{Z2act}
J_{j', h}(\tau, z).
\end{eqnarray}
This says that the $r$-tuple $(J_{1, h}, \hdots, J_{r, h})$ is a \emph{vector-valued
weak Jacobi form} of weight $0$ and index $1/2\langle h, h \rangle$. (By
(\ref{lambdacong}) we have $\langle h, h \rangle \in \mathbb{Q}$.) Part (i), which we do not need,
is proved by making use of a theorem of Miyamoto \cite{M}, which itself extends
some of the ideas in Zhu's modular-invariance theorem \cite{Z}. The proof of (ii) involves  applications of the ideas of Section 9, and in particular the permutation in (\ref{Z2act}) is the same as the one we described earlier (loc. cit.)

\pagebreak
\noindent
11. {\bf The quadratic space $(E, \langle \ , \ \rangle)$}

\medskip
We will prove the following result.
\begin{eqnarray}\label{Equad}
&&(E, \langle \ , \ \rangle)\ \mbox{is a \emph{positive-definite} rational quadratic space} \notag\\
&&\mbox{of rank $l$, and $L_0 \subseteq E$ is an additive subgroup of rank $l$.}
\end{eqnarray}

We have seen that both $E/L_0$ and $L_0/L$ are torsion groups. Hence
$E/L$ is also a torsion group, so in proving that $\langle h, h \rangle >0$ for $0 \not= h \in E$,
 it suffices to prove this under the additional assumption that
 $h \in L$. We assume this from now on, and set $m = \langle h, h \rangle$.
  Note that the results of Section 10 apply in this situation.
 
 \medskip
 We will show that $m \leq 0$ leads to a contradiction.
From (\ref{Z2act}) we know that $J_{j, h}(\tau, z+u\tau+v) = J_{j', h}(\tau, z)$
for $1 \leq j \leq r.$ In terms of the Fourier series $J_{j, h}:= \sum_{n, t} c(n, t)q^n\zeta^t, 
J_{j', h}:= \sum_{n, t} c'(n, t)q^n\zeta^t$, this reads
\begin{eqnarray*}
q^{\lambda_j-c/24}\sum_{n\geq 0, t} c(n, t)q^{n+mu^2/2+tu}\zeta^{t+mu} = q^{\lambda_{j'}-c/24}
\sum_{n \geq 0, t} c'(n, t)q^n\zeta^t 
\end{eqnarray*}
for all $u \in \mathbb{Z}$, $j'$ depending on $u$. Suppose first that $m=0$. If for some $t\not= 0$
there is $c(n, t) \not= 0$ we let $u \rightarrow -\infty$ and obtain a contradiction.
Therefore, $c(n, t)=0$ whenever $t\not= 0$. This says precisely that
$h_j(0)$ is the zero operator on $M^j$. Furthermore, this argument holds for any index $j$. But now
(\ref{Zmodinv2}) is contradicted. If $m<0$ the argument is even easier since we just have to let 
$u \rightarrow -\infty$ to get a contradiction.

\medskip
This proves that $\langle \ , \ \rangle$ is positive-definite on $E$, while rationality has already been established (\ref{lambdacong}). Now we prove that $E$ has rank $l$, using an argument familiar from the theory of root systems (cf.\ \cite{H}, Section 8.5).
 We have already seen (cf.\ (\ref{Etor}) and the line following) that $E$ contains a basis of $H$, say
$\{\alpha_1, \hdots, \alpha_l\}$.
We assert that $\{\alpha_1, \hdots, \alpha_l\}$ is a $\mathbb{Q}$-basis of $E$.

\medskip
Let $u \in E$. There are scalars
$c_1, \hdots, c_l \in\mathbb{C}$ such that $u = \sum_j c_j\alpha_j$.  We have for $1 \leq i \leq l$ that
\begin{eqnarray}\label{sole1}
\langle u, \alpha_i \rangle = \sum_j c_j \langle \alpha_i, \alpha_j \rangle.
\end{eqnarray}
Each $\langle u, \alpha_i \rangle$ and $\langle \alpha_i, \alpha_j \rangle$ are rational, and
the nondegeneracy of $\langle \ , \ \rangle$ implies that $(\langle \alpha_i, \alpha_j \rangle)$
is \emph{nonsingular}. Therefore,
$c_j = \langle u, \alpha_j \rangle/\det(\langle \alpha_i, \alpha_j \rangle) \in \mathbb{Q}$,
as required.

\medskip
We have proved that $E$ is a $\mathbb{Q}$-form for $H$, i.e.\ $H = \mathbb{C}\otimes_{\mathbb{Q}} E$, so that $E$ indeed has rank $l$. That $L_0 \subseteq E$ is a lattice of the same rank follows  from
(\ref{Etor}). All parts of (\ref{Equad}) are now established.

\medskip
Now  observe that the analysis that leads to the proof of (\ref{Equad}) carries over
 \emph{verbatim} to any \emph{nondegenerate} subspace $U \subseteq H$, say of rank $l'$. For such a subspace
we set $E' := U \cap E, L':= U\cap L, L_0':= U \cap L_0$. The result can then be stated as follows:
\begin{eqnarray}\label{Uquad}
&&(U, \langle \ , \ \rangle)\ \mbox{is a \emph{positive-definite} rational quadratic space} \notag\\
&&\mbox{of rank $l'$, and $L'_0 \subseteq E$ is an additive subgroup of rank $l'$.}
\end{eqnarray}

\medskip
Another application of weak Jacobi forms allows us to usefully strengthen the statement (\ref{Zmodinv2}) in some cases:
\begin{eqnarray}\label{morenontrivact}
\mbox{if $0\not= h \in E$ then}\ h_j(0)\not= 0\ \mbox{for each}\ 1 \leq j \leq r.
\end{eqnarray}
Suppose false. Because $E/L$ is a torsion group there is $0 \not= h \in L$ with $h_j(0)=0$ for some index $j$. Let $m = \langle h, h \rangle$, so that $m \not= 0$. Then $J_{j, h}(\tau, z)$ is a
pure $q$-expansion, i.e.\ no nonzero powers of $\zeta$ occur in the Fourier expansion. Indeed, it is just the partition function for $M^j$, so it also does not vanish.
By (\ref{Z2act}),
$e^{-\pi i \langle h, h \rangle(u^2 \tau +2uz)} J_{j', h}(\tau, z) = q^{-mu^2/2}\zeta^{-mu}\sum_{n\geq 0, t}
c'(n, t)q^{n-\lambda_j}\zeta^t$ is also a pure $q$-expansion. (As usual, $j'$ depends on $u$.)
 But because $m >0$ we can let $u \rightarrow \infty$ to see that in fact this power series is \emph{not} a pure $q$-expansion. This contradiction proves (\ref{morenontrivact}).

\bigskip \noindent
12. {\bf Commutants}

\medskip
We retain previous notation. In particular, from now on we fix a Cartan subalgebra $H \subseteq V_1$
and a nondegenerate subspace $U \subseteq H$ of rank $l'$. 
Let $M_U=(\langle U \rangle, Y, \mathbf{1}, \omega_U)$ be the Heisenberg
subVOA of rank $l'$ generated by $U$ (cf.\ (\ref{moreomegadefs})).
We set 
$Y(\omega_U, z) := \sum_{n \in \mathbb{Z}} L_U(n)z^{-n-2}$.

\medskip
Consider
\begin{eqnarray}\label{PUdef}
\mathcal{P}_U := \{(A, Y, \mathbf{1}, \omega_U) \ |\  A \subseteq V\}.
\end{eqnarray}
In words, $\mathcal{P}_U$ is the set of subVOAs $A \subseteq V$ which 
 have conformal vector $\omega_U$. $\mathcal{P}_U$ is partially ordered by
inclusion. It contains $M_U$, for example.

\medskip
One easily checks that $L(1)\omega_U=0$. Therefore, the theory of commutants
(\cite{FZ}, \cite{LL}, Section 3.11) shows that each $A \in \mathcal{P}_U$ has a 
\emph{compatible grading} with (\ref{specdecomp1}). That is 
$A_n := \{v \in A \ | \ L_U(0)v=n\} = A \cap V_n$. Moreover, $\mathcal{P}_U$ has a
\emph{unique maximal element}. Indeed, the \emph{commutant} $C_V(A) =$ ker$_VL_U(-1)$
for $A \in \mathcal{P}_U$ is \emph{independent} of $A$, and the maximal element of $\mathcal{P}_U$
is the double commutant $C_V(C_V(A))$.

\bigskip
\noindent
13. {\bf $U$-weights}

\medskip
Thanks to (\ref{crH}) we can use the language of weights to describe the action
of $u(0)\ (u \in U)$.  For $\beta \in U$ set
\begin{eqnarray*}
V(\beta) := \{ w \in V \ | \ u(0)w = \langle \beta, u\rangle w\ (u \in U)\}.
\end{eqnarray*}
$\beta$ is a \emph{$U$-weight},
or simply \emph{weight} (of $V$) if $V(\beta) \not= 0$, $V(\beta)$ is the $\beta$-weight space, and a nonzero $w \in V(\beta)$ is a weight vector of weight $\beta$.

\medskip
Using the action of $Y(u, z)\ (u \in U)$ on weight spaces, one
 shows that the set of $U$-weights
 \begin{eqnarray}\label{wtdef}
P := \{ \beta \in U \ |\ V(\beta) \not= 0\}
\end{eqnarray}
is a \emph{subgroup} of $U$. See \cite{DM1}, Section 4 for further details.
By the complete reducibility of $u(0)\ (u \in U)$
and the Stone von-Neumann theorem (\cite{FLM}, Section 1.7) applied to the Heisenberg
subVOA $M_U$, there is a weight space decomposition 
\begin{eqnarray}\label{Wspacedecomp}
V &=& M_U \otimes \Omega  = \bigoplus_{\beta \in P} M_U \otimes \Omega(\beta)
\end{eqnarray}
where $\Omega := \{v \in V \ | \ u(n)v=0\ (u \in U, n \geq 1)\}, \Omega(\beta) := \Omega \cap V(\beta)$, and $V(\beta) = M_U \otimes \Omega(\beta)$.

\medskip
$\Omega(0)$ is the commutant $C_V(M_U)$, and $M_U \otimes \Omega(0)$ the zero weight
space. By arguments in \cite{DM4} one sees that 
 $\Omega(0)$ is \emph{simple} VOA (the simplicity of $M_U$ is well-known), moreover each $V(\beta)$ is an \emph{irreducible}
 $M_U \otimes \Omega(0)$-module. So there is a tensor decomposition 
 \begin{eqnarray*}
V(\beta) = M_U(\beta) \otimes \Omega(\beta)
\end{eqnarray*}
where $M_U(\beta), \Omega(\beta)$ are irreducible modules  for $M_U, \Omega(0)$ respectively.
Furthermore, $V(\beta) \cong V(\beta')$ if, and only if, $\beta = \beta'$. In particular, there is an identification
\begin{eqnarray}\label{ebetadef}
M_U(\beta)= M_U \otimes e^{\beta}
\end{eqnarray}
where $e^{\beta} \in \Omega(\beta)$.

\bigskip
\noindent
14. {\bf Lattice subalgebras of $V$}

\medskip
We keep previous notation. In particular, $P$ is the group of $U$-weights (\ref{wtdef}) and
$E' = U \cap E, L' = L \cap U, L'_0 = L_0 \cap U$ are as in Section 11 (cf.\ (\ref{Uquad})).

\medskip
Since $(E', \langle \ , \ \rangle)$ is a rational space (\ref{Uquad}) and contains a basis of $U$, it follows that
$E' = \{u \in U \ | \ \langle u, E' \rangle \subseteq \mathbb{Q}\}.$ Because
$E'/L_0'$ is a torsion group, we then see that $(L'_0)^0 \subseteq E'.$ (Here, and below,
we set  $F^0:= \{ u \in U\ | \ \langle u, F \rangle \subseteq \mathbb{Z}\}$ for $F \subseteq E'$.)
Now $u \in P^0 \Leftrightarrow \langle P, u \rangle \subseteq \mathbb{Z} \Leftrightarrow$ all eigenvalues of $u(0)$ are integral $\Leftrightarrow u \in L'_0$. We conclude that
\begin{eqnarray}\label{Locontain}
P = (L'_0)^0 \subseteq E'
\end{eqnarray}
We will establish
\begin{eqnarray}\label{VLsubalgebras}
&&\mbox{there is a positive-definite even lattice $\Lambda \subseteq P$ such that $|P: \Lambda|$} \notag\\
&&\mbox{is finite and the maximal element $W$ of $\mathcal{P}_U$ satisfies $W \cong V_{\Lambda}$.} 
\end{eqnarray}
The argument utilizes ideas in \cite{DM2}.
Recall the isomorphism (\ref{Liiso}), which holds for all $u \in L_0$. Set
\begin{eqnarray}\label{moreLambdadef}
\Gamma:= \{ u \in L_0' \ | \ (V, Y_{\Delta_u(z)}) \cong (V, Y)\}.
\end{eqnarray}
This is a subgroup of $L_0'$ of finite index. Although not necessary at this stage, we can show immediately that 
$\Gamma$ is an even lattice. Indeed, if $u \in \Gamma$ then
the proof of (\ref{uurat1}) shows that we have $\lambda_j=0$ in that display, whence
$\langle u, u \rangle = n_0$ is a (nonnegative) integer. Now the assertion about $\Gamma$ follows from (\ref{Equad}).

\medskip
 There is another approach that gives more information. The isomorphism of $V$-modules
defined for $u \in \Gamma$ by (\ref{moreLambdadef}) implies the following assertion concerning the weight spaces
 in (\ref{Wspacedecomp}):
\begin{eqnarray}\label{Omegaisos}
\Omega(\beta) \cong \Omega(\beta+u)\ \ (u \in \Gamma, \beta \in P).
\end{eqnarray}
 In particular, taking $\beta=0$ shows that $\Omega(u) \not= 0\ (u \in \Gamma)$, whence
$\Gamma \subseteq P$. Using (\ref{Locontain}) we deduce
\begin{eqnarray*}
\Gamma  \subseteq P, P^0 \subseteq \Gamma^0,
\end{eqnarray*}
so $\Gamma$ is necessarily a positive-definite integral lattice of rank $l'$, and 
$|P: \Gamma| =:d$ is finite.
 (\ref{Omegaisos}) leads to a refinement of (\ref{Wspacedecomp}), namely a decomposition
 of $V$ into simple
 $M_U \otimes \Omega(0)$-modules
\begin{eqnarray*}
V = \bigoplus_{i=1}^d \bigoplus_{\beta \in \Gamma} M_U(\beta+\gamma_i)\otimes \Omega(\gamma_i),
\end{eqnarray*}
where $\{\gamma_i \ | \ 1 \leq i \leq d\}$ are coset representatives for 
$P/\Gamma$.

\medskip
Let
\begin{eqnarray}
&&\  \Lambda:= \{\beta \in P \ | \ \Omega(\beta)=\Omega(0)\}, \label{Lamdef} \\
&&W: = \bigoplus_{\beta \in \Lambda} M_U(\beta). \label{moreWdef}
\end{eqnarray}
Then $\Gamma \subseteq \Lambda$ and $W = C_V(\Omega(0)) = C_V(C_V(M_U))$. In particular,
$\Lambda$ is an additive subgroup of $P$ of finite index and $W$ is a subVOA of $V$. Indeed, it is the maximal element of the poset $\mathcal{P}_U$ discussed in Section 12.

\medskip
The $L(0)$-weight of $e^{\beta} \in W (\beta \in \Lambda)$ (cf.\ (\ref{ebetadef})) coincides with its
$L_U(0)$-weight (cf.\ Section 12). Using the associativity formula, we have
\begin{eqnarray*}
L(0)e^{\beta}&=&L_0(0).e^{\beta} = 1/2\sum_{t=1}^{l'} (h_t(-1)h_t)(1)e^{\beta} \\
&=&1/2\sum_{t=1}^{l'} \left\{\sum_{k \geq 0} h_t(-1-k)h_t(1+k)+h_t(-k)h_t(k) \right\}e^{\beta} \\
&=&1/2\sum_{t=1}^{l'} h_t(0)h_t(0)e^{\beta} = 1/2\sum_{t=1}^{l'} \langle \beta, h_t \rangle^2 e^{\beta} = 1/2\langle\beta, \beta\rangle,
\end{eqnarray*}
showing that $1/2\langle \beta, \beta \rangle \in \mathbb{Z}\ (\beta \in \Lambda)$. 

\medskip
This shows that $\Lambda$ is an \emph{even} lattice of rank $l'$.
 The isomorphism $W \cong V_{\Lambda}$
then follows from the uniqueness of simple current extensions (\cite{DM1}, Section 5). 
 This completes the proof of (\ref{VLsubalgebras}) and Theorem 1 is established.

\bigskip
\noindent
15. {\bf The tripartite subVOA of $V$}

\medskip
We consider more closely the consequences of (\ref{VLsubalgebras}) in the case that $U:=$ rad$(V_1)$ is the radical of $V_1$. (\ref{VLsubalgebras}) is applicable here because $A$ is indeed nondegenerate (cf.\ Section 4).  We keep the notation from previous Sections.

\medskip
Observe that in this case, the lattice $\Lambda$ \emph{contains no roots}, i.e.\ there is no
$\beta \in \Lambda$ satisfying $\langle \beta, \beta \rangle =2$. For if $\beta \in \Lambda$ is a root 
then $\beta$ is contained in an $sl_2$-subalgebra of $V_1$ and hence cannot lie in $A$.
The commutant $\Omega(0)$ of $W$ contains the 
Levi factor $S \subseteq V_1$, hence also the subVOA $G$ that it generates (cf.\ (\ref{Wiso})).
We can then consider the commutant of $G$ in $\Omega(0)$, call it $Z$. In this way we obtain
the canonical conformal subVOA of $V$ that we call the \emph{tripartite subalgebra}
\begin{eqnarray*}
T = W \otimes G \otimes Z.
\end{eqnarray*}

 \medskip
By construction, $T$ is a \emph{conformal subalgebra} of $V$, and 
the conformal gradings on $W, G, Z$ are compatible with the $L(0)$-grading on $V$.
Because $(W \otimes G)_1 = W_1 \oplus G_1 = V_1$ then $Z_1=0$.  This completes
the proof of Theorem 3.

\bigskip
\noindent
{\bf Conjecture}: $Z$ is a strongly regular VOA.

\medskip
This is just a special case of a more general conjecture, namely that the commutant of a rational subVOA (in a strongly regular VOA, say) is itself rational. If the Conjecture is true then the tripartite subalgebra $T$ is strongly regular, and $V$ reduces to a \emph{finite sum} of irreducible $T$-modules. In this way, the classification of strongly regular VOAs reduces to the classification of strongly regular VOAs
$Z$ with $Z_1=0$ and the extension problem as discussed in the Introduction.

\bigskip \noindent
16. {\bf The invariants $\tilde{c}$ and $l$}

\medskip
We give some applications of  Theorem 1 exemplifying the philosophy of the
previous paragraph. Let $V$ be a strongly regular VOA of central charge $c$
and $H\subseteq V$ a 
Cartan subalgebra of rank $l$. Recall (\cite{DM1}) that the \emph{effective central charge}
of $V$ is the quantity
\begin{eqnarray*}
\tilde{c}_V=\tilde{c} := c - 24\lambda_{min}.
\end{eqnarray*}
Here, $\lambda_{min}$ is  the \emph{minimum} of the conformal weights
$\lambda_j\ (1 \leq j \leq r)$ of the irreducible $V$-modules. It is known (loc.\ cit.) that $\tilde{c}\geq l$
and  $\tilde{c}>0$ if $\dim V>1$. Because of these facts,
$\tilde{c}$ is often a more useful invariant  than $c$ itself. Note that $\tilde{c}$ is defined
for any rational VOA. 

\medskip
We now give the proof of Theorem 7. 
The basic idea is to combine Zhu's modular-invariance
\cite{Z} together with growth conditions on the Fourier coefficients of components of vector-valued modular forms
\cite{KnM}. This method was first used in \cite{DM1}. The availability of Theorem 1 brings added clarity.

\medskip
It follows easily from the definitions that if
$W \subseteq V$ is a \emph{conformal} subalgebra then $\tilde{c}_V \leq \tilde{c}_W$.
Moreover $\tilde{c}$ is multiplicative over tensor products (\cite{FHL}, Section 4.6).
So if (b) of Theorem 7 holds then  $\tilde{c}_V \leq \tilde{c}_{V_{\Lambda}}+ \tilde{c}_{L(c_{p, q}, 0)}$.
Since rk$\Lambda = l$ then $\tilde{c}_{V_{\Lambda}} = c = l$ because $\lambda_{min}=0$
for lattice theories \cite{D}. Moreover, for the discrete series
Virasoro VOA we have
 (\cite{DM1}, Section 4, Example (e))
\begin{eqnarray}\label{tildecval}
\tilde{c} = 1 - \frac{6}{pq} \ \ ((p, q)=1, 2\leq p < q),
\end{eqnarray}
in particular we always have $\tilde{c}_{L(c_{p, q}, 0)}<1$.
 Therefore $\tilde{c}_V < l+1$. This establishes the implication (b) $\Rightarrow$ (a) in Theorem 7.

\medskip
Next, taking $U=H$  in
Theorem 1, we find that the maximal element of $\mathcal{P}_H$ is a lattice subVOA
$W\cong V_{\Lambda}$ with rk$\Lambda = \dim H = l$. Let $C=C_V(W)$ be the commutant
of $W$. Then $W \otimes C$ is a conformal subVOA of $V$. Now suppose that part 
(b) of the Theorem does \emph{not} hold. Thus the Virasoro subalgebra of $C$, call it $Vir_C$, has a central charge $c'$, say, that is \emph{not} 
 in the discrete series.
Then the known submodule structure of
Verma
modules over the Virasoro algebra shows that the \emph{partition function}
$Z_{Vir_C}(\tau)$ of $Vir_C$ satisfies
\begin{eqnarray*}
Z_{Vir_C}(\tau) := \mbox{Tr}_{Vir_C}q^{L(0)-c'/24} = q^{-c'/24}\prod_{n=2}^{\infty}(1-q^n)^{-1}.
\end{eqnarray*}
(\cite{DM1}, Proposition 6.1 summarizes exactly what we need here.) Therefore,
\begin{eqnarray*}
Z_{W\otimes Vir_C}(\tau) &:=& Z_W(\tau) Z_{Vir_C}(\tau) \\
&=&
\frac{\theta_{\Lambda}(\tau)}{\eta(\tau)^l} \frac{q^{-c'/24}}{\prod_{n=2}^{\infty}(1-q^n)} \\
&=&
\frac{\theta_{\Lambda}(\tau)}{\eta(\tau)^{l+1}} q^{(1-c')/24}(1-q).
\end{eqnarray*}
($\theta_{\Lambda}(\tau)$ and $\eta(\tau)$ are the \emph{theta-function} of $\Lambda$ and the 
\emph{eta-function} respectively.) It follows that for any $\epsilon>0$, the coefficients of the $q$-expansion of 
$\eta(\tau)^{l+1-\epsilon}Z_{W\otimes Vir_C}(\tau)$ have \emph{exponential growth}. Therefore, the same statement holds true \emph{ipso facto} if we replace $W \otimes Vir_C$ with $V$.
We state this as
\begin{eqnarray}\label{expgrowth}
\mbox{the coefficients of $\eta(\tau)^{l+1-\epsilon}Z_V(\tau)$ have \emph{exponential growth} 
($\epsilon>0$).}
\end{eqnarray}
 
On the other hand, consider the column vector 
\begin{eqnarray*}
F(\tau):= (Z_{M^1}(\tau), \hdots, Z_{M^r}(\tau))^t
\end{eqnarray*}
whose components are the partition functions of the irreducible
$V$-modules $M^j$. By Zhu's modular-invariance theorem, $F(\tau)$ is a vector-valued modular form
of weight $0$ on the full modular group $SL_2(\mathbb{Z})$ associated with some representation
of $SL_2(\mathbb{Z})$. (See \cite{MT}, Section 8 for a discussion of vector-valued modular forms
in the context of VOAs.) Moreover, each $Z_{M^j}(\tau)$ is holomorphic
in the complex upper half-plane, so that their only poles are at the cusps. The very definition of
$\tilde{c}$, and the reason for its importance,
 is that the maximum order of a pole of \emph{any} of the partition functions $Z_{M^j}(\tau)$
 is $\tilde{c}/24$. It follows from this that
 \begin{eqnarray*}
\eta(\tau)^{\tilde{c}}F(\tau)
\end{eqnarray*}
is a \emph{holomorphic} vector-valued modular form on $SL_2(\mathbb{Z})$.
As such, the Fourier coefficients of the component functions have \emph{polynomial growth}
\cite{KnM}.
In particular, this applies to $\eta(\tau)^{\tilde{c}}Z_V(\tau)$, which is one of the components.

\medskip
Comparing the last statement with (\ref{expgrowth}), it follows that $\tilde{c} > l+1-\epsilon$
for all $\epsilon>0$, i.e.\ $\tilde{c} \geq l+1$. So we have shown that if  part (b) of the Theorem does not hold, neither does part (a). Theorem 7 is thus proved.

\medskip
The special case $l=0$ of the Theorem characterizes \emph{minimal models}. 
We state it as 

\medskip \noindent
{\bf Theorem 8}: Let $V$ be a strongly regular VOA. Then $\tilde{c}<1$ if,
and only if, the Virasoro subalgebra of $V$ is in the discrete series. 

\medskip \noindent
{\bf Corollary 9}: Let $V$ be a strongly regular VOA with $\dim V > 1$. Then
$\tilde{c}\geq 2/5$, and equality holds if, and only if, $V \cong L(c_{2, 5}, 0)$,
the (Yang-Lee) discrete series Virasoro VOA with $c=-22/5$.

\medskip
Because $\dim V >1$ then $\tilde{c}>0$, and if $\tilde{c}<1$
then  $V$ is a minimal model by Theorem 8.
Inspection of (\ref{tildecval}) shows that the least positive value is
$2/5$, corresponding to the Yang-Lee model. This theory has only two
irreducible modules, of conformal weight $0$ and $-1/5$. Therefore
the second irreducible cannot be contained in $V$, so that
$V \cong L(c_{2, 5}, 0)$, as asserted in Corollary 9. Informally, the Corollary says that
the Yang-Lee theory is the \emph{smallest} rational CFT.

\medskip
We give a final numerical example. Suppose that $V$ is a strongly regular simple VOA such that $1<\tilde{c}<7/5$.
Since $l\leq c$ we must have $l=0$ or $1$. In the latter case, by Theorem 7 we see that $V$ contains as a conformal subVOA a tensor product $V_{\Lambda} \otimes Vir$ where $Vir$ is a Virasoro algebra in the discrete series with $0<\tilde{c}_{Vir}<2/5$. This is impossible by Corollary 9. So in fact $l=0$,
 i.e.\ $V$ has Lie rank $0$, meaning that $V_1=0$. The smallest value of $\tilde {c}$ in the range
 $(1, 7/5)$  that I know of is a parafermion theory with $\tilde{c} = 8/7$.

\end{document}